\documentclass[a4paper,12pt]{article}
\usepackage{amssymb}
\usepackage{amsthm}
\usepackage{amsxtra}
\usepackage{amsfonts}
\usepackage[all]{xy}
\usepackage[T1]{fontenc}
\newcommand{\PRCB}{T}

\newcommand{\LOn}{\mathcal{L}_n^*}
\newcommand{\LC}{\mathcal{L}_n}

\newcommand{\Li}{\mathcal{L}^{(i)}}

\newcommand{\clos}[2]{#1^{(#2)}}

\newcommand{\tp}{\mathrm{tp}}
\newcommand{\qftp}{\mathrm{qftp}}

\newtheorem{thm}{Theorem}[section]

\newtheorem{lem}[thm]{Lemma}

\newtheorem{prop}[thm]{Proposition}

\newtheorem{fact}[thm]{Fact}

\theoremstyle{definition}
\newtheorem{para}[thm]{}
\newtheorem{Types:}[thm]{Types:}

\newtheorem*{claim}{Claim}
\newtheorem*{claim1}{Claim 1}
\newtheorem*{claim2}{Claim 2}

\newtheorem{defn}[thm]{Definition}

\newtheorem{rem}[thm]{Remark}

\newtheorem{notconv}[thm]{Notations and Conventions}

\def\rest{{\lower      .25    em      \hbox{$\vert$}}}

\font\helpp=cmsy5
\def\semdp
{\hbox{$\times\kern-.23em\lower-.1em\hbox{\helpp\char'152}$}\,}

\def\dnfo{\,\raise.2em\hbox{$\,\mathrel|\kern-.9em\lower.35em\hbox{$\smile$}
$}}
\def\dnf#1{\lower1em\hbox{$\buildrel\dnfo\over{\scriptstyle #1}$}}
\def\dfo{\;\raise.2em\hbox{$\mathrel|\kern-.9em\lower.35em\hbox{$\smile$}
\kern-.7em\hbox{\char'57}$}\;}
\def\df#1{\lower1em\hbox{$\buildrel\dfo\over{\scriptstyle #1}$}}

\def\Ind#1#2{#1\setbox0=\hbox{$#1x$}\kern\wd0\hbox to 0pt{\hss$#1\mid$\hss}     
\lower.9\ht0\hbox to 0pt{\hss$#1\smile$\hss}\kern\wd0}
\def\Notind#1#2{#1\setbox0=\hbox{$#1x$}\kern\wd0\hbox to 0pt{\mathchardef
\nn=12854\hss$#1\nn$\kern1.4\wd0\hss}\hbox to
0pt{\hss$#1\mid$\hss}\lower.9\ht0 \hbox to
0pt{\hss$#1\smile$\hss}\kern\wd0}

\def\tho{\text{\th}}                                        

\begin{document}

\title{Imaginaries in bounded pseudo real closed fields}
\author{Samaria Montenegro}
\date{}
\maketitle

\begin{abstract}
The main result of this paper is that if $M$ is a bounded $PRC$ field then $Th(M)$ eliminates imaginaries in the language of rings expanded by constant symbols.

\end{abstract}

\section{Introduction}
A  \emph{pseudo algebraically closed field} ($PAC$ field) is a field $M$ such that every absolutely irreducible affine variety defined over $M$ has an $M$-rational point.
The concept of a $PAC$ field was introduced by J.Ax in \cite{Ax} and has been extensively studied.
The above definition of $PAC$ field has an equivalent model-theoretic version: $M$ is existentially closed (in the language of rings) into each regular field extension of $M$. 

A field $M$ is called \emph{bounded} if for any integer $n$, $M$ has only finitely many extensions of degree $n$.
Hrushovski showed in \cite{Hrus} that if $M$ is a bounded $PAC$ field, and $\mathcal{L}$ is the language of rings expanded by enough constants, then $Th_{\mathcal{L}}(M)$ eliminate imaginaries.

The notion of $PAC$ field has been generalized  by Basarab  in  \cite{Ba0} and then by Prestel in \cite{Pre} for ordered fields. 
Prestel calls a field $M$ \emph{pseudo real closed field} ($PRC$ field) if $M$ is existentially closed (in the language of rings) into each regular field extension $L$ to which all orderings of $M$ extend. 
Remark that if $M$ is a $PRC$ field and has no orderings, then $M$ is a $PAC$ field.
$PRC$ fields were extensively studied by L. van den Dries in \cite{Van}, A. Prestel in \cite{Pre}, M. Jarden in \cite{J1}, \cite{J2}, \cite{J3}, S. Basarab in \cite{Ba} and \cite{Ba1}, and others.

The main result in this paper is a generalization to $PRC$ bounded fields of elimination of imaginaries for $PAC$ fields.

As corollary of the elimination of imaginaries and the fact that the algebraic closure in the sense of model theory defines a pregeometry we obtain (Theorem \ref{PRCrosy}) that the complete theory of a bounded $PRC$ field is superrosy of $U^{\tho}$ -rank $1$.

The organization of the paper is as follows:
In section \ref{PRCpreliminaires} we give the required preliminaries on pseudo real closed fields and we fix a complete theory $\PRCB$ of a bounded $PRC$ field, where we enrich the language adding constants for an elementary submodel. 
In section \ref{VOn} for $n \geq 1$, we define the theory $VO_n$ in a multi-sorted language $\mathcal{L}_n^*$.
To each model of the theory $\PRCB$ we associate a model of $VO_n$ (Remark \ref{PRCOVn}).
This result is an important tool in the proof of elimination of imaginaries for bounded $PRC$ fields.
We show quantifier elimination and elimination of imaginaries for the theory $VO_n$. 
Finally in section \ref{EIsecPRC} we prove the elimination of imaginaries for bounded $PRC$ fields (Theorem \ref{EIPRC}).

\bigskip
\noindent \textbf{Acknowledgments}
I would like to thank Zo\'e Chatzidakis for her support and guidance through this work, and for all her corrections and suggestions.
\section{Pseudo real closed fields}\label{PRCpreliminaires}

In this section we give the required preliminaries in pseudo real closed fields.

\begin{notconv}
If $M$ is a model of an $\mathcal{L}$-theory $T$ and $A \subseteq M$, then $\mathcal{L}(A)$ denotes the set of $\mathcal{L}$-formulas with parameters in $A$. 
If $\bar{a}$ is a tuple of $M$, we denote by $\tp_{\mathcal{L}}^M(a/A)$ $(\qftp_{\mathcal{L}}^M(a/A))$ the set of $\mathcal{L}(A)$-formulas (quantifier-free $\mathcal{L}(A)$-formulas) $\varphi$, such that $M \models \varphi(a)$. 
Denote by $acl^M_{\mathcal{L}}$ and $dcl^M_\mathcal{L}$ the model theoretic algebraic and definable closure in $M$. 
We omit $M$ or $\mathcal{L}$ when the structure or the language is clear.
We denote by $\mathcal{L}_{\mathcal{R}}$ the language of rings.
\bf{All fields considered in this paper will have characteristic zero}.

\end{notconv}

\begin{defn}Let $M,N$ be fields such that $M \subseteq N$.
\begin{enumerate}
 \item The extension $N/M$ is called \emph{totally real} if each order on $M$ extends to some order on $N$.
 \item We say that $N$ is a \emph{regular extension} of $M$ if $N \cap M^{alg}=M$.
\end{enumerate}

\end{defn}

\begin{fact} \label{PRC} \cite[Theorem 1.2]{Pre}
For a field $M$ the following are equivalent:
\begin{enumerate}
	\item $M$ is existentially closed (relative to $\mathcal{L_R}$) in every totally real regular extension $N$ of $M$.
	\item For every absolutely irreducible variety $V$ defined over $M$, if $V$ has a simple 
	 $\overline{M}^r-$rational point for every real closure $\overline{M}^r$ of $M$, then $V$ has an $M$-rational point.
\end{enumerate}
\end{fact}

\begin{defn} 
A field $M$ that satisfies the conditions of Fact \ref{PRC} is \emph{pseudo real closed} ($PRC$). 
By \cite[Theorem 4.1]{Pre} we can axiomatize the class of $PRC$ fields in  $\mathcal{L_R}$. 
Remark that the class of $PRC$ fields contains the class of $PAC$ fields and the class of real closed fields ($RCF$ fields).

In the case when $M$ admits only a finite number of orderings this already implies that $M$ is existentially closed in $N$ even in the language augmented by predicates for each order $<$ of $M$ \cite[Theorem 1.7]{Pre}.

\end{defn}

\begin{fact} \cite[Proposition 1.2]{J2}\label{PRCcaracte}
Let $M$ be a $PRC$ field. Then:
\begin{enumerate}
	\item If $<$ is an order on $M$, then $M$ is $<$-dense in $\overline{M}^r$, where $\overline{M}^r$ is the real closure of $M$ with respect to the order $<$. 
	\item If $<_i$ and $<_j$ are different orders on $M$, then $<_i$ and $<_j$ induce different topologies.
	\end{enumerate}
\end{fact}

\begin{para}\label{ApTh} 
\textbf{Approximation Theorem} \cite[1.7]{Van} 
Let $M$ be a field  and $\tau_1, \ldots, \tau_n$ different topologies on $M$ induced by orders. 
For each $i \in \{1,..,n\}$, let $U_i$  be a non-empty $\tau_i$-open subset of $M$.
Then $\displaystyle{\bigcap_{i=1}^n U_i \not = \emptyset}$.
\end{para}

\subsection{Bounded pseudo real closed fields}

\begin{para}
\textbf{Bounded fields:} 
If $M$ is a bounded field and $M^*$ an elementary extension of $M$ (in a language containing the language of rings) then the restriction map: $G(M^*)\rightarrow G(M)$ is an isomorphism \cite[Lemma 1.22]{Cha0}. 
\end{para}

\begin{para}\label{PRCB}
\textbf{Notation:}
We fix a bounded $PRC$ field $K$, which is not real closed and a countable elementary substructure  $K_0$ of $K$. 
Then $G(K_0)\cong G(K)$ and $K_0^{alg}K= K^{alg}$.
Since $K$ is bounded there exists $n \in \mathbb{N}$ such that $K$ has exactly $n$ distinct orders (see Remark 3.2 of \cite{Mon}).
Then we can view $K$ as a structure of the form $(K, <_1, \ldots, <_n)$, where $\{<_1, \ldots, <_n\}$ are all the different orders on the field $K$.

In this paper we will work over $K_0$, thus we denote by $\mathcal{L}$ the language of rings with constant symbols for the elements of $K_0$, $\Li:= \mathcal{L} \cup \{<_i\}$ and $\mathcal{L}_n:= \mathcal{L} \cup \{<_1, \ldots, <_n\}$.
We let $\PRCB:=Th_{\LC}(K)$.

If $n = 0$, then $K$ is a $PAC$ field and by Corollary 3.1 of \cite{Hrus} $Th_{\mathcal{L}}(M)$ has elimination of imaginaries.
Thus we will suppose that $n\geq 1$.

\end{para}

\begin{fact}\label{PRCacl}\cite{Mon}
Let $(M, <_1, \ldots, <_n)$ be a model of $\PRCB$.
\begin{enumerate}
\item For all $i \in \{1, \ldots,n\}$ the order $<_i$ is $\exists$-definable in the language $\mathcal{L}$, and $Th_{\mathcal{L}}(M)$ is model complete.
\item If  $A$ is a subfield of $M$ containing $K_0$, then $A^{alg}\cap M = acl_{\LC}^{M}(A)= dcl_{\LC}^M(A) = acl_{\mathcal{L}}^{M}(A)= dcl_{\mathcal{L}}^M(A)$.
 \end{enumerate}

\end{fact}

\section{The theory $V\!O_n$}\label{VOn}

\begin{defn}
Fix $n \in \mathbb{N}$, $n \geq 1$.
Let $\LOn$ be the $(n+1)$-sorted language consisting of $n+1$ sorts $\{R_0, \ldots, R_{n}\}$, $n$ binary relation symbols $\{<_1, \ldots, <_n\}$, with $<_i$ on the sort $R_i$, a constant symbol $0$ of sort $R_0$, and $2n$ function symbols $\{f_1, \ldots, f_n, g_1, \ldots, g_n\}$, where $f_i: R_0 \rightarrow R_i$ and $g_i: R_i \rightarrow R_0$.
Let $V\!O_n$ be the $\LOn$-theory axiomatized by:
\begin{enumerate}
\item $<_i$ defines a dense linear order without endpoints on $R_i$, for $i \in \{1, \ldots,n\}$,
\item $f_i: R_0 \rightarrow R_i$ is an injective function, for $i \in \{1, \ldots,n\}$,
\item $g_i(x) := \left\{ \begin{array}{ll}
         f_i^{-1}(x) & \mbox{if $x \in f_i(R_0)$};\\
        0 & \mbox{otherwise}.\end{array} \right. $
\item $f_i(R_0)$ is $<_i$-dense and $<_i$-co-dense in $R_i$, for $i \in \{1, \ldots,n\}$,
\item If $a_i,b_i \in R_i$ are such that $a_i <_i b_i$, then there exists $x_0 \in R_0$ such that $a_i <_i f_i(x_0) <_i b_i$, for all $i \in \{1, \ldots,n\}$.

\end{enumerate}

\textbf{Notation}: \begin{enumerate}
\item  If $a_i,b_i \in R_i \cup \{\pm \infty\}$ and $a_i<_ib_i$, denote by $(a_i,b_i)_i:= \{x \in R_i: a_i<_ix<_ib_i\}$ and by $f_i^{-1}(a_i,b_i)_i:=\{x \in R_0: f_i(x)\in (a_i,b_i)_i\}$.
\item If $M$ is a model of $VO_n$ and $A \subseteq M$, denote by $\langle A \rangle$ the $\LOn$-substructure generated by $A$.
\end{enumerate}

\begin{rem}\label{PRCOVn}
If $M$ is a model of $T$ (see \ref{PRCB}), we associate to $M$ a model of $V\!O_n$ as follows:
For each $i \in \{1, \ldots, n\}$ we let $\clos{M}{i}$ be a real closure of $M$ for the order $<_i$.
We define $R_0=M$, $R_i=\clos{M}{i}$ with its natural order $<_i$, $f_i:M \rightarrow \clos{M}{i}$ the natural inclusion and $g_i$ its ``inverse''.
Axioms $1,2$ and $3$ are clearly true. 
Axiom $5$ is true by the Approximation Theorem (\ref{ApTh}) and Fact \ref{PRCcaracte}. 
For axiom $4$, by Fact \ref{PRCcaracte}(1), $M$ is $<_i$-dense in $\clos{M}{i}$. 
To see that $M$ is co-dense in $\clos{M}{i}$, suppose by contradiction that there exists an $<_i$-interval $I$ such that $I\subseteq M$.
Let $\alpha \in \clos{M}{i} \setminus M$ and define  the $<_i$-interval $J:= \{x + \alpha: x \in I\}$.
Then $J\subset (\clos{M}{i}\setminus M)$, which contradicts $M$ being dense in $\clos{M}{i}$.

Observe that if $a_i,b_i \in \clos{M}{i} \cup \{\pm \infty\}$ and $I= \displaystyle\bigcap_{i=1}^n((a_i,b_i)_i\cap M)$, then in the corresponding model of $V\!O_n$, the set $I$ is $\LOn$-definable by the formula $\displaystyle{\bigcap_{i=1}^nf_i^{-1}(a_i,b_i)_i}$.  
\end{rem}

\begin{thm}\label{EQVOn}
The theory $V\!O_n$ is $\aleph_0$-categorical, has quantifier elimination in $\LOn$ and if $M$ is a model of $V\!O_n$ and $A \subseteq M$, then   $acl_{\LOn}(A)=dcl_{\LOn}(A)=\langle A \rangle= \displaystyle{\bigcup_{a \in A} \langle a \rangle }$.
\begin{proof}
Let $M=(M_0,\ldots, M_n)$ and  $N=(N_0,\ldots, N_n)$ be countable models of $V\!O_n$. 
Let $\mathcal{I}(M,N)$ be the family of partial isomorphisms with finite domain of $M$ to $N$.
Let $g \in \mathcal{I}(M,N)$, $A = dom(g)$ and $B= g(A)$. 
Let $b \in M \setminus A$; we need to find $c \in M$ such that if $\tilde{g}: A \cup \{b\}\rightarrow B \cup \{c\}$ is defined by $\tilde{g}|_{A}= g$ and $\tilde{g}(b)=c$, then $\tilde{g} \in \mathcal{I}(M,N)$. 

As $g$ extends uniquely to an isomorphism from $\langle A \rangle$ to $\langle B \rangle$, we can suppose that $A =\langle A \rangle$ and $B=\langle B \rangle$.

\textbf{Case 1:} $b \in M_0.$

Observe that $\langle A,b \rangle = A  \cup \{b, f_i(b): 1 \leq i \le n\}$.
For each $i \in \{1, \ldots,n\}$, let  $I^i$ be the smallest open $<_i$-interval containing $f_i(b)$ with extremities in $(A \cap M_i) \cup \{\pm \infty\}$.
Let $J^i\subseteq N_i$ be the interval obtained by applying $g$ to the extremities of $I^i$. 
By axiom 5 there exists $c \in N_0$ such that  $f_i(c) \in J^i$ for all $1 \leq i \leq n$; define $g(b)=c$.

\textbf{Case 2:} There exists $i>0$ such that $b \in M_i.$ 

There exists a smallest open $<_i$-interval $I^i$  with extremities in $(A \cap M_i) \cup \{\pm \infty\}$, such that $b \in I^i$. 
Let $J^i\subseteq N_i$ be the interval obtained by applying $g$ to the extremities of $I^i$. 
We have two possibilities:
\begin{enumerate}
 \item $b \in f_i(M_0)$: apply Case 1 to $g_i^{-1}(b)=b'$, observe that $\langle A,b \rangle = A \cup \{b', f_j(b'): 1 \leq j \leq n\}$.
 \item $b \not\in f_i(M_0)$: as $f_i(N_0)$ is $<_i$- co-dense in $N_i$, there is $c \in J^i \setminus f_i(N_0)$. Define $g(b)=c$. Then $\langle A,b \rangle = A \cup \{b\}$.
\end{enumerate}
This shows quantifier elimination and $\aleph_0$-categoricity. The assertion about the algebraic and definable closures is clear.
\end{proof}
\end{thm}

\end{defn}

\begin{defn} Let $(M_0, \ldots, M_n)$ be a model of $V\!O_n$.
\begin{enumerate}
\item A set of the form $I= \displaystyle\bigcap_{i=1}^n{f_i^{-1}}(I^i)$, with $I^i$ a non-empty $<_i$-open interval in $M_i$ is called a \emph{multi-interval}. 
 Observe that if a multi-interval $I= \displaystyle\bigcap_{i=1}^n{f_i^{-1}}(I^i)$ is definable over $A$, then by quantifier elimination (Theorem \ref{EQVOn}) each $I^i$ has its extremities in $(\langle A \rangle \cap M_i) \cup \{\pm \infty\}$. 
 We call the \emph{set of extremities of the multi-interval $I$} the set of extremities of $I^i$, for all $i\in \{1, \ldots,n\}$.
 \item Let $E \subseteq M_0$. 
 We say that $E$ is \emph{multi-open} if for each $e \in E$, there exists a multi-interval $I$ such that $e \in I$ and $I  \subseteq E$.
\item  Let $E \subseteq M_0$ and $e \in E$. 
 Let $I= \displaystyle{\bigcap_{i=1}^nf_i^{-1}(I^i)}$ be a multi-interval such that $e \in I $ and $I \subseteq E$.
We say that $I$ is the \emph{maximal multi-interval in $E$ containing $e$} if for all $m \in  \{1, \ldots n\}$, $I^m$ is the maximal $<_m$-interval with the property that:
\begin{enumerate}
 \item $f_m(e) \in I^m$,
 \item There are $J^{m+1}, \ldots, J^{n}$ intervals in $M_{m+1}, \ldots, M_{n}$ respectively, such that:
 \begin{enumerate}
 \item $f_j(e) \in J^j$, for all $m+1 \leq j \leq n$,
  \item $\displaystyle \bigcap_{l=1}^m f_l^{-1}({I^l}) \cap \bigcap_{j=m+1}^n f_j^{-1}(J^j) \subseteq E$.
 \end{enumerate}
\end{enumerate}

\end{enumerate}
\end{defn}

\begin{rem}\label{unicitycd}
Observe that if $X(e)= \displaystyle{\bigcap_{i=1}^nf_i^{-1}(X^i(e))}$ and $Y(e)= \displaystyle{\bigcap_{i=1}^nf_i^{-1}(Y^i(e))}$ are maximal multi-intervals in $E$ containing $e$,
then for all $i \in \{1, \ldots, n\}$, $X^i(e)= Y^i(e)$: 
It is clear using maximality that $X^1(e)= Y^1(e)$; by induction suppose that $X^1(e)=Y^1(e), \ldots, X^m(e)=Y^m(e)$, using maximality again we obtain that $X^{m+1}(e)= Y^{m+1}(e)$.
\end{rem}

\begin{lem}\label{maxinter}
Let $M=(M_0, \ldots, M_n)$ be a model of $V\!O_n$ and $A \subseteq M$ finite.
Let $E \subseteq M_0$ be multi-open $\LOn(A)$-definable, and let $e \in E$. 
Then there exists a maximal multi-interval $X(e)$ in $E$ containing $e$ and its extremities are in $\langle A \rangle \cup \{\pm \infty\}$.
\begin{proof}
Let $Y^1(e)$ be the set of open $<_1$-intervals $J^1$ satisfying:
\begin{enumerate}
 \item $f_1(e) \in J^1.$
 \item There exist $J^2, \ldots, J^n$ open intervals in $M_2, \ldots, M_n$ respectively such that:
 \begin{enumerate}
  \item $f_j(e) \in J^j$, for all $2\leq j \leq n$,
  \item $\displaystyle{\bigcap_{j=1}^nf_j^{-1}(J^j)} \subseteq E.$
 \end{enumerate}

\end{enumerate}

Observe that if $J^1, L^1 \in Y^1(e)$ then $J^1 \cup L^1 \in Y^1(e)$.

\begin{claim} {$Y^1(e)$ has a maximal element.}
\begin{proof}
Since $e \in E$ and $E$ is multi-open, using axiom $5$ we can find $e_1, e_2 \in M_0$, such that for all $i \leq n$, $f_i(e_1)<_if_i(e_2)$, and $e \in \displaystyle{\bigcap_{i=1}^n}f_i^{-1}(f_i(e_1),f_i(e_2))_i \subseteq E.$
Then $(f_1(e_1),f_1(e_2))_1 \in Y^1(e)$.
Define $X^1(e,e_1,e_2):= \{x \in M_1: \exists a,b \in M_1 (x \in (a,b)_1 \wedge f_1(e) \in (a,b)_1 \wedge f_1^{-1}(a,b)_1 \cap \displaystyle{\bigcap_{i=2}^nf_i^{-1}(f_i(e_1),f_i(e_2))_i} \subseteq E \}$. 

Observe that $X^1(e,e_1,e_2)$ is definable, $(f_1(e_1),f_1(e_2))_1 \subseteq X^1(e,e_1,e_2)$ and that it is connected for the $<_1$-topology: if $x,y \in X^1(e,e_1,e_2)$ and $x <_1 y$ then $(x,y)_1 \subseteq X^1(e,e_1,e_2)$. 
This implies by quantifier elimination (Theorem \ref{EQVOn}) that $X^1(e,e_1,e_2)$ is an $<_1$-interval, so that $X^1(e,e_1,e_2) \in Y^1(e)$.

Since $X^1(e,e_1,e_2)$ is definable with parameters in $A \cup \{e,e_1,e_2\}$, its extremities are in $dcl_{\LOn}(A \cup \{e,e_1,e_2\}) \cup \{\pm \infty\}$.
So by Theorem \ref{EQVOn}, its extremities are in $(\langle A \rangle \cap M_1) \cup \{f_1(e), f_1(e_1),f_1(e_2)\} \cup \{\pm\infty\}.$ 
As $f_1(e)\in (f_1(e_1),f_1(e_2))_1 \subseteq X^1(e,e_1,e_2)$, it cannot be one of the extremities. 
Thus $X^1(e,e_1,e_2)$ has its extremities in  $(\langle A \rangle \cap M_1)  \cup \{f_1(e_1),f_1(e_2)\} \cup \{\pm\infty\}$.

Let $e'_1, e'_2 \in M_0$ be such that $(f_i(e'_1),f_i(e'_2))_i \subsetneqq (f_i(e_1),f_i(e_2))_i$ for all $i \leq n$, and $e \in \displaystyle{\bigcap_{i=1}^n}f_i^{-1}(f_i(e'_1),f_i(e'_2))_i \subseteq E.$
Then $X^1(e,e_1,e_2)\subseteq X^1(e,e_1',e_2')$.

As before we obtain that $X^1(e,e_1',e_2')$ is an $<_1$-open interval with extremities in $(\langle A \rangle \cap M_1)  \cup \{f_1(e'_1),f_1(e'_2)\} \cup \{\pm\infty\}$ and that $X^1(e,e_1',e_2') \in Y^1(e)$.
But as $f_1(e'_1),f_1(e'_2) \in (f_1(e_1),f_1(e_2))_1 \subseteq X^1(e,e_1,e_2)\subseteq X^1(e,e_1',e_2')$, then $f_1(e'_1),f_1(e'_2)$ cannot be the extremities of $X^1(e,e_1',e_2')$. 
This implies that $X^1(e,e_1',e_2')$ has its extremities in $(\langle A \rangle \cap M_1) \cup \pm \{\infty\}$.

We have shown that if $e'_1, e'_2 \in M_0$ are such that $(f_i(e'_1),f_i(e'_2))_i \subsetneqq (f_i(e_1),f_i(e_2))_i$, then $X(e,e_1',e_2')$ is $\LOn(A)$-definable and that as the multi-intervals $(e_1',e_2')_i$ decrease, the sets $X^1(e,e_1,e_2)$ increase.
By $\aleph_0$-categoricity and the fact that $\langle A \rangle $ is finite, the sets $X^1(e,e_1',e_2')$ stabilize for some value of $e_1',e_2'$ and so $Y^1(e)$ has a maximal element $X^1(e)$, which has its extremities in $\langle A \rangle \cup \{\pm \infty\}$.

\end{proof}
\end{claim}
By induction suppose that we have already defined $Y^1(e), \ldots, Y^{l-1}(e)$, for $l \leq n$, with $X^1(e), \ldots, X^{l-1}(e)$ the maximal element of $Y^1(e), \ldots, Y^{l-1}(e)$ respectively.

Let $Y^l(e)$ be the set of $<_l$-intervals $J^l$ satisfying:
\begin{enumerate}
 \item $f_l(e) \in J^l$,
 \item There exist $J^{l+1}, \ldots, J^n$ intervals in $M_{l+1}, \ldots, M_n$ respectively such that:
 \begin{enumerate}
  \item $f_j(e) \in J^j$, for all $l+1\leq j \leq n$,
  \item  $\displaystyle{\bigcap_{j=1}^{l-1}f_j^{-1}(X^{j}(e)) \cap \bigcap_{j=l}^nf_j^{-1}(J^j)} \subseteq E$.
 \end{enumerate}
\end{enumerate}
Reasoning exactly as for $l=1$, and defining 
\[\begin{split}
X^l(e,e_1,e_2):= \{x \in M_l: \exists a,b \in M_l (x \in (a,b)_l \wedge  f_l(e) \in (a,b)_l \wedge\\ \displaystyle{\bigcap_{i=1}^{l-1}X^i(e)}\cap (a,b)_l  \cap \displaystyle{\bigcap_{i=l+1}^nf_i^{-1}(f_i(e_1),f_i(e_2))_i} \subseteq E \},
\end{split}\]
we find that $Y^l(e)$ has a maximal element $X^l(e)$, which is $\LOn(A)$-definable.  

Let $X(e):= \displaystyle{\bigcap_{i=1}^nf_i^{-1}(X^i(e))}$. Then $X(e)$ is the maximal multi-interval in $E$ containing $e$, and it is $\LOn(A)$-definable. 
\end{proof}
\end{lem}

\begin{para} \textbf{Canonical decomposition:}

Let $M= (M_0, \ldots,M_n)$ be a model of $V\!O_n$ and $E\subseteq M_0$. 
Let $B \subseteq \mathbb{N}$ and $E_0 \subseteq E$ be finite sets and for all $j \in B$, let $I_j$ be a multi-interval.
We say that $\displaystyle{\bigcup_{j\in B}I_j} \cup E_0$ is a \emph{canonical decomposition of $E$} if: 
\begin{enumerate}
 \item $E = \displaystyle{\bigcup_{j\in B}I_j} \cup E_0$,
 \item for all $e \in E_0$, there is no multi-interval $I$ containing $e$ such that $I \subseteq E$,
 \item for all $e \in E\setminus E_0$, there exists $j \in B$ such that $I_j$ is the maximal multi-interval in $E$ containing $e$,
 \item for all $j \in B$ there exists $e \in  E \setminus E_0$ such that $I_j$ is the maximal multi-interval in $E$ containing $e$.
\end{enumerate}
 
\end{para}

\begin{thm}\label{candescom}
Let $M= (M_0, \ldots,M_n)$ be a model of $V\!O_n$ and $A \subseteq M$ finite. Let $E \subseteq M_0$ be $\LOn(A)$-definable. 
Then there exists a unique canonical decomposition of $E$ and its extremities are in $\langle A \rangle \cup \{\pm \infty\}$.
\begin{proof}
Define $\widetilde{E}:= \{x \in E: \mbox{there exists a multi-interval} \; I\; \mbox{such that}\; x \in I\; \mbox{and}\; I \subseteq E \}$ and $E_0:= E \setminus{\widetilde{E}}$.
Observe that $\widetilde{E}$ and $E_0$ are $\LOn(A)$-definable, that $\widetilde{E}$ is multi-open and that $E = \widetilde{E} \cup E_0$. 
Using quantifier elimination (Theorem \ref{EQVOn}) and the fact that the negation of an atomic formula is a disjunction of atomic formulas we obtain that $E_0$ is a finite set, defined by disjunctions of equalities.

As $\widetilde{E}$ is multi-open, by Lemma \ref{maxinter} for each $e \in \widetilde{E}$ there exists $X(e)= \displaystyle{\bigcap_{i=1}^n f_i^{-1}(X^i(e))}$, the maximal multi-interval in $\widetilde{E}$ containing $e$, and its extremities are in $\langle A \rangle \cup \{\pm \infty\}$.

Since $\langle A \rangle$ is finite then $\{X(e): e \in \widetilde{E}\}$ is finite.
Let $B \subset \widetilde{E}$ be finite such that $\{X(e): e \in \widetilde{E}\}= \{X(e): e \in B\}$.
Then $\widetilde{E} = \displaystyle{\bigcup_{e\in B}X(e)}$.

Therefore $\displaystyle{\bigcup_{e\in B}X(e)} \cup E_0$ is a canonical decomposition of $E$ and its extremities are in $\langle A \rangle \cup \{\pm \infty\}$.
The uniqueness is clear by Remark \ref{unicitycd}.
\end{proof}
\end{thm}

\begin{rem}\label{remdescomcan}
 The uniqueness in Theorem \ref{candescom} implies that if $E \subseteq M_0$ is definable with parameters in $\langle A \rangle$ and also with parameters in $\langle B \rangle$, then the canonical decomposition of $E$ is definable with parameters in $\langle A \rangle \cap \langle B \rangle$.
 
 \end{rem}

\begin{thm}
 $V\!O_n$ has elimination of imaginaries in the language $\LOn$.
 \begin{proof}
Let $M=(M_0, \ldots, M_n)$ be a model of $V\!O_n$.
 \begin{claim}
$V\!O_n$ has unary elimination of imaginaries: 
\begin{proof}
Let $A \subseteq M$ be finite and let $E$ be an $\LOn(A)$-definable set.
We have two cases:

Case 1: $E \subseteq M_0$.

By Theorem \ref{candescom} there exists a canonical decomposition of $E$. 
Let $\bar{c}$ be the set of finite points and extremities of the multi-intervals in the canonical decomposition of $E$. 
By Remark \ref{remdescomcan}, $\bar{c}$ is the code of the set $E$.

Case 2: $E \subseteq M_i$ for some $i>0$.

Define $E_1:= E \cap f_i(M_0)$, $E_2:= E \cap (M_i \setminus f_i(M_0))$.
Observe that $E_1$ and $E_2$ are $\LOn(A)$-definable and that $E$ is the disjoint union of $E_1$ and $E_2$.
Since $g_i(E_1)\subseteq M_0$, by Case 1 $g_i(E_1)$ is coded by some tuple $c_1 \in \langle A \rangle$.

By quantifier elimination (Theorem \ref{EQVOn}), $E_2$ is defined by a boolean combination of formulas of the form: $a<_i x$, $a= x$, $b <_j f_jg_i(x)$, $b = f_jg_i(x)$ with $a \in M_i \cap \langle A \rangle$, $b \in M_j \cap \langle A \rangle$. 
Since $E_2 \subseteq M_i \setminus f_i(M_0)$, then $f_j(g_i(x))=f_j(0)$ for all $x \in E_2$. 
This implies that $E_2$ is defined by a formula $\psi(x) \wedge g_i(x)= 0$, where $\psi(x)$ is a finite union of points and disjoint intervals with extremities in $(\langle A \rangle \cap M_i) \cup \{\pm \infty\}$. 
Then the set $c_2$ of  finite points and extremities of these intervals is the code of the set $E_2$, and $c:= (c_1,c_2)$ is the code of $E$.
\end{proof}
\end{claim}

By Remark $3.2.2$ of \cite{HHM} it is enough to show that every definable unary function (with parameters) is encoded in $M$.
Let $i,j \in \{0, \ldots,n\}$ and $h: M_i \rightarrow M_j$ be an $\LOn(A)$-definable function with $A \subseteq M$ finite.

Let $B:= \{x: h(x) \in \langle x \rangle \}$.
Observe that $h(M_i \setminus B)$ is finite:
Let $x \in M_i \setminus B$; since $h$ is $\LOn(A)$-definable, then $h(x) \in dcl(Ax) = \langle A \rangle \cup  \langle x \rangle$. 
 Since $x \not \in B$ then $h(x) \not \in  \langle x \rangle$, and $h(x) \in \langle A \rangle \cap M_j$. 
As $\langle A \rangle$ is finite, $h(M_i \setminus B) \subseteq \langle A \rangle$ is finite.

Let $m \in \mathbb{N}$ and $\{a_1, \ldots, a_m\} \in \langle A \rangle \cap M_j$ be such that $h(M_i \setminus B)= \{a_1, \ldots, a_m\}$.
Let $X_l:= \{x \in M_i: h(x)=a_l\}$; by unary elimination the set $X_l$ is coded in $M$ by a tuple $\ulcorner X_l \urcorner$.

The function $h|_{B}$ is also coded in $M$:
Let $X_0:= \{x \in M_i: h(x) \in \langle 0 \rangle\} \subseteq B$; by unary elimination the set $X_0$ is coded in $M$ by a tuple $\ulcorner X_0 \urcorner$.
If $D:= B \setminus X_0$, then
\[h|_D = \left\{ \begin{array}{ll}
         id & \mbox{if}\; i=j=0,\\
        f_j g_i & \mbox{if}\; i,j>0,\\
        f_j & \mbox{if}\; i=0, j>0,\\
        g_j & \mbox{if}\; i>0, j=0.\\
        \end{array} \right. \]

Let $d:= (\ulcorner X_0 \urcorner, \ulcorner X_1 \urcorner, \ldots, \ulcorner X_m \urcorner, a_1, \ldots, a_m)\subseteq M$.
Then $d$ is the code in $M$ of $h$.
\end{proof}
  \end{thm}

\section{Elimination of Imaginaries in bounded PRC fields}\label{EIsecPRC}

\begin{defn}
Let $(M, <_1, \ldots, <_n)$ be a model of $\PRCB$ (see \ref{PRCB}).
Denote by $\clos{M}{i}$ a fixed real closure of $M$ with respect to $<_i$.
\begin{enumerate}
\item A subset of $M$ of the form $I= \displaystyle{\bigcap_{i=1}^n (I^i\cap M)}$ with $I^i$ a non-empty $<_i$-open interval in $\clos{M}{i}$ is called a \emph{multi-interval}.
Observe that by  \ref{ApTh} (Approximation Theorem) and Fact \ref{PRCcaracte} (1) every multi-interval is non empty.

\item A definable subset $S$ of $M$ is called \emph{multi-open} if for each $x \in S$, there exist a multi-interval $I$ such that $x \in I$ and $I \subseteq M$.

\item A definable subset $S$ of a multi-interval $I = \displaystyle{\bigcap_{i=1}^n (I^i\cap M)}$ is called \emph{multi-dense} in $I$ if for any multi-interval $J \subseteq I$, $J \cap S \not = \emptyset.$
Note that multi-density implies $<_i$- density in $I^i$, for all $i \in \{1. \ldots,n\}$.
\end{enumerate}
\end{defn}

\begin{rem}\label{definM}
Let $(M, <_1, \ldots, <_n)$ be a model of \PRCB. 
Let $i \in \{1, \ldots,n\}$ and $a \in \clos{M}{i} \setminus M$ such that $a \in acl^{\clos{M}{i}}(c)$, with $c$ a tuple in $M$. 
Then $A= \{x \in M: x <_i a\}$ is definable  in $M$ by a quantifier-free $\Li(c)$-formula. 
\begin{proof}
By quantifier elimination of the theory of real closed fields (RCF) and the fact that $acl^{\clos{M}{i}}= dcl^{\clos{M}{i}}$, we can find a quantifier-free $\Li$-formula $\phi(x,c)$, such that $\clos{M}{i}\models \forall x (x <_i a \leftrightarrow \phi(x,c))$.
Then $x \in A$ if and only if $M\models \phi(x,c)$.
 \end{proof}
\end{rem}

\begin{fact} \label{descomposition}\cite[Theorem 3.13]{Mon}
Let $(M, <_1, \ldots, <_n)$ be a model of $\PRCB$ and $A \subseteq M$.  Let $S \subseteq M$ be an $\LC(A)$-definable set.
Then there are a finite set $S_0\subseteq S$, $m \in \mathbb{N}$ and $I_1, \ldots, I_m$, with $I_j = \displaystyle{\bigcap_{i=1}^n (I^i_{j}\cap M)}$ a multi-interval for all $j \in \{1, \ldots,m\}$ such that:
\begin{enumerate}
\item $S \subseteq \displaystyle{\bigcup_{j=1}^m I_ j \cup S_0},$
\item$\{x\in I_j: x \in S\}$ is multi-dense in $I_{j}$ for all  $1 \leq j \leq  m$,
\item $I^i_{j} \subseteq \clos{M}{i}$ has its extremities in $dcl^{\clos{M}{i}}_{\Li}(A)\cup \{\pm \infty\}$ for all  $1 \leq j \leq  m$  and $1 \leq i \leq n$, 
\item the set $I^i_{j}\cap M$ is definable in $M$ by a quantifier-free $\Li(A)$-formula, for all $1 \leq j \leq  m$ and $1 \leq i \leq n$.
 \end{enumerate}
\end{fact}

\begin{prop}\label{descompositioncan} 
Let $(M, <_1, \ldots, <_n)$ be a model of $T$ and $A,B \subseteq M$. 
Let $S \subseteq M$, $\LC$-definable with parameters in $A$ and also in $B$. 
Then there are a finite set $S_0 \subseteq S$, $m\in \mathbb{N}$ and $I_1, \ldots, I_m$, with $I_j:= \displaystyle{\bigcap_{i=1}^n}(I^i_j\cap M)$ a multi-interval such that:
\begin{enumerate}
 \item $S \subseteq \displaystyle{\bigcup_{j=1}^m I_j} \cup S_0$,
\item $\{x \in I_j: x \in S\}$ is multi-dense in $I_j$, for all $1 \leq j \leq m$,
\item the set $I_j^i \cap M$ is definable in $M$ by a quantifier-free $\Li$-formula with parameters in $acl(A) \cap acl(B)$. 
 \end{enumerate}
 \begin{proof}
 Define $\widetilde{S}:= \{x \in M: \mbox{there exists a multi-interval}\; I\; \mbox{such that} \; x \in I \; \mbox{and}\; I\cap S \; \mbox{is multi-dense in} \;I \}$.
 Since $S$ is $\LC(A)$-definable and also $\LC(B)$-definable, then $\widetilde{S}$ is $\LC(A)$-definable and also $\LC(B)$-definable.

By Fact \ref{descomposition} there exists a finite set $S_0 \subseteq S$ such that $S \subseteq \widetilde{S} \cup S_0$. 
As $\widetilde{S}$ is multi-open and $\LC(A)$-definable, using Fact \ref{descomposition} there exists $r_1 \in \mathbb{N}$ and multi-intervals $I_j= \displaystyle{\bigcap_{i=1}^n(I_j^i\cap M)}$, for all $1 \leq j \leq r_1$ such that:
$\widetilde{S}= \displaystyle{\bigcup_{j=1}^{r_1}I_j}$, and for all $1 \leq i \leq n$, $1 \leq j \leq r_1$, $I^i_j$ has its extremities in $dcl_{\Li}^{\clos{M}{i}}(A) \cup \{\pm \infty\}$. 
 
Similarly, as $\widetilde{S}$ is also $\LC(B)$-definable, there exists $r_2 \in \mathbb{N}$ and multi-intervals $J_j= \displaystyle{\bigcap (J_j^i\cap M)}$ for all $1 \leq j \leq r_2$ such that:
$\widetilde{S}= \displaystyle{\bigcup_{j=1}^{r_2}J_j}$, and for all $1 \leq i \leq n$, $1 \leq j \leq r_2$, $J^i_j$  has its extremities in $dcl_{\Li}^{\clos{M}{i}}(B) \cup \{\pm \infty\}$. 
  
Let $\widetilde{A}$ be the set of extremities of $I^i_j$, for $1 \leq i \leq n$ and $1 \leq j \leq r_1$, and let $\widetilde{B}$ be the set of extremities of $J^i_j$, for $1 \leq i \leq n$ and $1 \leq j \leq r_2$.
 
Consider the $\LOn$-structure $\widetilde{M}=(M, \clos{M}{1}, \ldots, \clos{M}{n})$ associated to $M$ (see Remark \ref{PRCOVn}).
Observe that $\widetilde{S}$ is $\LOn(\widetilde{A})$-definable and also $\LOn(\widetilde{B})$-definable in $\widetilde{M}$.
Then by Theorem \ref{candescom} and Remark \ref{remdescomcan} there exists a unique canonical decomposition of $\widetilde{S}$ and it is $\LOn$-definable with parameters in $\langle \widetilde{A}\rangle \cap \langle \widetilde{B}\rangle$. 

Let $m \in \mathbb{N}$ and $I_1, \ldots I_m$ the multi-intervals such that $\displaystyle{\bigcup_{j=1}^m I_j}$ is the canonical decomposition of $\widetilde{S}$.
If $I_j= \displaystyle \bigcap_{i=1}^n (I_j^i\cap M)$, by Remark \ref{definM} and the definition of $\widetilde{A}$ and $\widetilde{B}$, $I_j^i \cap M$ is definable in $M$ by a quantifier-free $\Li$-formula with parameters in $dcl_{\LC}(A) \cap dcl_{\LC}(B)$.
 
 \end{proof}
\end{prop}

\begin{lem}\label{indimag}
Let $M$ be a model of $T$ and $e \in M^{eq}$. Let $a$ be a tuple in $M$ and $f$ an $\LC(\emptyset)$-definable function such that $f(a)=e$. 
Let $E = acl(E) \supseteq  acl^{eq}(e)\cap M$.
Then there exist tuples $b,b'$ in $M$, $ACF$-independent over $E$, such that $\tp(b/E)=  \tp(b'/E) = \tp(a/E)$ and $f(a)=f(b)= f(b')$.
\begin{proof}
 The proof is exactly the same as in claim 1 of Proposition 3.1 in \cite{Hrus}.
 If $a \in E$, it is clear. Assume $a \not \in E$.
 
 Sketch: Remember that by Fact \ref{PRCacl} if $A \subseteq M$, then $acl(A)= dcl(A)= A^{alg} \cap M$.
 Using Neumann's Lemma we can find conjugates $a_1,a_2$ of $a$ over $E \cup \{e\}$ satisfying:
 \[acl(E,a_1) \cap acl(E,a_2)= E.\]
Take such $a_1, a_2$ with $trdeg(a_2/Ea_1)=m$ maximal satisfying \eqref{eq1} below\\ 
\begin{equation} \label{eq1}
\tp(a_1/E)= \tp(a_2/E),\; acl(E,a_1) \cap acl(E,a_2)= E,\; \mbox{and} \; f(a_1)=f(a_2)=e. \tag{1}
\end{equation}

Take $a_3$ $ACF$-independent of $a_2$ over $E(a_1)$, such that $\tp(a_3/Ea_1)= \tp(a_2/Ea_1)$.
Then $f(a_3)= f(a_2)= f(a_1)=e$ and $acl(E,a_1) \cap acl(E,a_3)= E$.

Since $a_3$ is $ACF$-independent of $a_2$ over $E(a_1)$, then $acl(E,a_3,a_1)\cap acl(E,a_2,a_1)=acl(E,a_1)$. 
Intersecting both sides with $acl(E,a_3)$ we obtain $acl(E,a_3)\cap acl(E,a_2,a_1) = E$ and then $acl(E,a_3)\cap acl(E,a_2)=E$.

Using the maximality of $m$, $trdeg(a_3/Ea_2) \leq m$, and so $a_3$ is $ACF$-independent of $a_2$ over $E(a_1)$ and it is also $ACF$-independent of $a_1$ over $E(a_2)$. 
Elimination of imaginaries in $ACF$ and the fact that $E(a_1)^{alg}\cap E(a_2)^{alg}= E^{alg}$ imply that $a_3$ is $ACF$-independent of $a_1a_2$ over $E$.
Let $b = a_1$ and $b'= a_3$. 
 \end{proof}
 
 \end{lem}

\begin{lem}\label{Imagdense}
Let $M$ be a sufficiently saturated model of $T$, $e \in M^{eq}$, $a \in M$ and $f$ an $\LC(\emptyset)$-definable function such that $f(a)=e$. 
Let $E = acl^{eq}(e)\cap M$. 
Suppose that $e \notin dcl^{eq}(E)$. 
Then there is a multi-interval $I = \displaystyle{\bigcap_{i=1}^n}{(I^i\cap M)}$ such that $a \in I$ and $\{x \in I: \tp(x/E) = \tp(a/E) \wedge f(x)\neq e\}$ is multi-dense in $I$.
\begin{proof}
By Lemma \ref{indimag} there exists $b \in M$, $ACF$-independent of $a$ over $E$, such that $\tp(b/E)= \tp(a/E)$ and $f(a)=f(b)$.
For each formula $\alpha(x) \in \tp(a/E)$, define $\Phi_{\alpha}(x,y):= \alpha(x) \wedge f(x) \not =f(y)$. 
Then $\Phi_{\alpha}(M,a)=\Phi_{\alpha}(M,b):= A_{\alpha}$. 

Since $e \notin dcl^{eq}(E)$, $\tp(a/E)$ is not algebraic and $\tp(a/E) \cup \{f(x)\neq e\}$ is consistent.
Take $d \in M$ realizing $\tp(a/E) \cup \{f(x)\neq e\}$.
Then $d \in  A_{\alpha}$.

Since $A_{\alpha}$ is $\LC(Ea)$-definable and also $\LC(Eb)$-definable, by Proposition \ref{descompositioncan} there exists a multi-interval $J_{\alpha}=\displaystyle\bigcap_{i=1}^n{(J^i_{\alpha}\cap M)}$ such that:
\begin{enumerate}
\item[1.]$d \in J_{\alpha}$
\item[2.]$\{x \in J_{\alpha}:  x \in A_{\alpha}\}$ is multi-dense in $J_{\alpha}$.
\item[3.] $J^i_{\alpha}\cap M$ is definable in $M$ by a quantifier-free $\Li$-formula with parameters in $acl(Ea)\cap acl(Eb)= E$.
\end{enumerate}

So $J_{\alpha}$ is $\LC(E)$-definable in $M$.
As $\tp(d/E)=tp(a/E)$ and $d \in J_{\alpha}$, then $a \in J_{\alpha}$.

By saturation and Fact \ref{PRCcaracte}(1), there exists for all $i \in \{1, \ldots,n\}$ an $<_i$-interval $I^i$, with extremities in $M$ such that $a \in I^i \subseteq \displaystyle{\bigcap_{\alpha(x) \in \tp(a/E)}(J^i_{\alpha}}\cap M)$. 
Then $ a \in I:= \displaystyle{\bigcap_{i=1}^nI^i}$ and $\{x \in I: M \models \alpha(x) \wedge f(x) \not = e\}$ is multi-dense in $I$, for all $\alpha(x) \in \tp(a/E)$.
This implies using saturation that $\{x \in I: \tp(x/E) = \tp(a/E) \wedge f(x)\neq e\}$ is multi-dense in $I$.
\end{proof}
\end{lem}

\begin{fact}\label{IT2}\cite[Theorem 3.21]{Mon}
Let $(M, <_1, \ldots, <_n)$ be a model of  $\PRCB$. Let $E = acl(E) \subseteq M$ and $a_1, a_2, d$ tuples of $M$ such that:
$d$ is ACF-independent of $\{a_1,a_2\}$ over $E$, $\tp_{\LC}(a_1/E)= \tp_{\LC}(a_2/E) $, and $\qftp_{\LC}(d,a_1/E) = \qftp_{\LC}(d,a_2/E) $. 
Suppose that $E(a_1)^{alg}\cap E(a_2)^{alg}= E^{alg}$.

Then there exists a tuple $d^{*}$ in some elementary extension $M^{*}$ of $M$ such that:
\begin{enumerate} 
\item $d^*$ is ACF-independent of $\{a_1,a_2\}$ over $E$,
\item $\tp_{\LC}(a_1, d^{*}/E) = \tp_{\LC}(a_2, d^{*}/E),$ 
\item $\tp_{\LC}(a_1, d^{*}/E)= \tp_{\LC}(a_1, d/E).$
\end{enumerate}
\end{fact}

\begin{thm}\label{EIPRC}
$T$ has elimination of imaginaries.
\begin{proof}
Since we are working with a field it is enough to show that $T$ has weak elimination of imaginaries.
Let $M$ be a monster model of $T$ and $e \in M^{eq}$. Define $E:= acl^{eq}(e)\cap M$. 
We need to show that $e \in dcl^{eq}(E)$. 
Let $a$ be a tuple from $M$ and let $f$ be an $\LC(\emptyset)$-definable function such that $f(a) = e$. 
Suppose that $e \not\in dcl^{eq}(E).$

\begin{claim1}We can suppose that $trdeg(E(a)/E)$=1:
\begin{proof}
Choose $a$ with $trdeg(E(a)/E)$ minimal such that $f(a)=e$. 
Take $a'\subseteq a$ such that $trdeg(E(a)/E(a'))=1.$

By Lemma \ref{indimag} there is a tuple $b$ in $M$, $ACF$-independent of $a$ over $E$, such that $\tp(a/E)=\tp(b/E)$ and $f(a)=f(b)$.

Since $b$ is $ACF$-independent of $a$ over $E$ and $a \notin acl(Ea')$, then $a \notin acl(Ea'b)$. As $e \in dcl^{eq}(b)$ then $acl^{eq}(Eea')\subseteq acl^{eq}(Ea'b)$. 
Thus $a \notin acl^{eq}(Eea')$. It follows without loss of generality that we can replace $E$ by $acl(E(a'))$.
\end{proof}
\end{claim1}

Suppose that $a = (a_1, \ldots, a_m)$, $a_1 \not \in E$ and $a \subseteq acl(Ea_1)$. 
Then $a_j \in acl(Ea_1)= dcl(Ea_1)$, for all $j \in \{1, \ldots, m\}$, and we can suppose that $m=1$.

By Lemma \ref{indimag} there exists $b \in M$, $ACF$-independent of $a$ over $E$, such that $\tp(a/E)=\tp(b/E)$ and $f(a)=f(b)$. 

By Lemma \ref{Imagdense} there is a multi-interval $I= \displaystyle{\bigcap_{i=1}^n (I^i\cap M)}$ such that $a \in I$, and $\{x \in I: \tp(x/E) = \tp(a/E) \wedge f(x)\neq e\}$ is multi-dense in $I$.

\begin{claim2}
$\qftp_{\LC}(a/Eb) \cup \tp_{\LC}(a/E) \cup \{f(x)\neq e\}$ is consistent:
\begin{proof}
By compactness it is enough to show that if $\psi(x,b)\in \qftp_{\LC}(a/Eb)$, then there exists $d$ such that $\tp_{\LC}(d/E)=\tp_{\LC}(a/E)$ and $M \models \psi(d,b) \wedge  f(d)\neq e$.
 
As $\psi(x,b)\in \qftp(a/Eb)$ and $a \notin acl(Eb)$, there is a multi-interval $J(b) := \displaystyle{\bigcap_{i=1}^n(J^i(b)\cap M)}$ such that $a \in J(b)$ and $J(b) \subseteq \psi(M,b)$.

As  $a \in I^i \cap J^i(b)$ we can assume by taking the intersection that $J^i(b)\subseteq I^i$, for all $i \in \{1, \ldots,n\}$. 
Since $\{x \in I: \tp(x/E) = \tp(a/E) \wedge f(x)\neq e\}$ is multi-dense in $I$ we can find $d \in J(b)$, such that $f(d)\neq e$, and $\tp(d/E) = \tp(a/E)$. 
Then $M \models \psi(d,b) \wedge f(d)\neq e$.
\end{proof}
\end{claim2}

Let $d$ realize $\qftp_{\LC}(a/Eb) \cup \tp_{\LC}(a/E) \cup \{f(x)\neq e\}$. Then $\qftp_{\LC}(b, d/E) = \qftp_{\LC}(b, a/E)$, and $f(d) \not= f(b)$. 
By Fact \ref{IT2} we can find $b^{*}$ such that $\tp_{\LC}(b^{*}, d/E)= \tp_{\LC}(b^{*}, a/E)$ and  $\tp_{\LC}(b^{*},a/E )= \tp_{\LC}(b, a/E)$. 
As $f(b)=f(a)=e$, then $f(b^{*})= f(a)= e$. But we have also that $f(b^{*})= f(d) \not= e$. This is a contradiction.
 
\end{proof}

\end{thm}

\begin{fact} \cite[Theorem 4.12]{EaOn} \label{EIrosy}
 Any theory which has geometric elimination of imaginaries and for which algebraic closure defines a pregeometry is superrosy and $U^{\tho}(x=x)= 1$.
\end{fact}

\begin{thm}\label{PRCrosy}
 Let $M$ be a $PRC$ bounded field and let $T= Th_{\mathcal{L}}(M)$ (see \ref{PRCB}). Then $T$ is superrosy and $U^{\tho}(x=x)= 1.$
 \begin{proof}
By Fact \ref{PRCacl}, algebraic closure defines a pregeometry and by Theorem \ref{EIPRC} $\PRCB$ has elimination of imaginaries. Then by Fact \ref{EIrosy} $\PRCB$ is superrosy and $U^{\tho}(x=x)= 1.$.
 \end{proof}

\end{thm}

\end{document}